\documentclass[11pt]{amsart}

\usepackage{amssymb,amsmath,amsthm,mathtools}
\usepackage[hmargin=28mm,vmargin=30mm]{geometry}
\usepackage{microtype,url}

\newcommand{\C}{\mathbb{C}}
\newcommand{\R}{\mathbb{R}}
\newcommand{\norm}[1]{\left\lVert #1\right\rVert}
\newcommand{\id}{\operatorname{Id}}
\renewcommand{\le}{\leqslant}
\renewcommand{\ge}{\geqslant}

\theoremstyle{plain}
\newtheorem{theorem}{Theorem}[section]
\newtheorem{lemma}[theorem]{Lemma}
\newtheorem{prop}[theorem]{Proposition}

\theoremstyle{remark}
\newtheorem{remark}[theorem]{Remark}

\numberwithin{equation}{section}

\begin{document}

\title[There is no universal separable Banach algebra]{There is no universal separable Banach algebra}
\author{Tomasz Kania}
\date{}
\dedicatory{Dedicated to W.~B.~Johnson on the occasion of his 80\textsuperscript{th} birthday}
\thanks{RVO: 67985840}

\begin{abstract}
We prove that no separable Banach algebra is universal for homomorphic embeddings of all separable Banach algebras, whether embeddings are merely bounded or required to be contractive. The same holds in the commutative category.

The proof is short and direct. For each closed subspace $X$ of $c_0$ we build a separable commutative Banach algebra $A_X$ whose multiplication records the canonical pairing between $X^*$ and $X$. Any injective homomorphism of $A_X$ into a Banach algebra $B$ forces the identity operator on $X^*$ to factor through $B$. Thus a universal separable Banach algebra would be complementably universal for the family of duals of subspaces of $c_0$, contrary to a theorem of Johnson and Szankowski.
\end{abstract}

\maketitle

\section{Introduction}

A natural question in the theory of Banach algebras concerns the existence of universal objects for certain classes thereof: does there exist a separable Banach algebra that contains, in a suitable sense, all separable Banach algebras from the given class? The answer may depend critically on both the algebraic structure considered and the notion of `containment' employed.

For commutative C$^*$-algebras, the situation is well understood. Every separable commutative C$^*$-algebra is of the form $C_0(K)$ for some locally compact metrisable space $K$. Choose a compact metrisable space $K^+$ containing $K$ as an open subset such that $K^+\setminus K=\{\infty\}$; when $K$ is non-compact, this is the one-point compactification, and when $K$ is compact, one may take the topological sum $K\sqcup\{\infty\}$. Since every compact metrisable space is a quotient of the Cantor set $\Delta := 2^{\mathbb N}$, there exists a continuous surjection $\Delta \to K^+$, and hence $C(K^+)$ embeds isometrically as a subalgebra of $C(\Delta)$. The algebra $C_0(K)$ identifies with the codimension-one ideal
\[
\{f\in C(K^+): f(\infty)=0\}
\]
of $C(K^+)$. Thus $C(\Delta)$ is universal for separable commutative C$^*$-algebras. In contrast, there is no separable C$^*$-algebra that is universal for all separable C$^*$-algebras, as can be deduced from \cite[Proposition~2.6]{JungePisier:1995}.

For general Banach algebras without additional structure, it is easy to see that there exists no separable commutative Banach algebra that is universal for isometric embeddings. Indeed, for each $t\ge 1$, let $A_t=\C^2$ with coordinatewise multiplication and norm
\[
\norm{(z,w)}_t := |z|+t|w|.
\]
Then $A_t$ is a commutative Banach algebra and $e_t:=(0,1)$ is an idempotent with $\norm{e_t}_t=t$. If a commutative Banach algebra $B$ contained isometric copies of all separable commutative Banach algebras, then for every $t\ge 1$ it would contain an idempotent $p_t$ with $\norm{p_t}=t$. For distinct $s,t\ge 1$, the idempotents $p_s$ and $p_t$ are distinct and commute, so $(p_s-p_t)^3=p_s-p_t$. Hence
\[
1\le \norm{p_s-p_t}^2,
\]
that is, $\norm{p_s-p_t}\ge 1$. Therefore $B$ would contain an uncountable $1$-separated set, impossible in a separable metric space.

However, to the best of our knowledge, the question of universality for bounded (or contractive) homomorphic embeddings, that is, embeddings required only to have closed range, has not been addressed in the literature. This is the problem studied in the present note.

We work over the field of complex numbers. Banach algebras are not assumed unital, and homomorphisms are not assumed unital. An \emph{embedding} is an injective homomorphism that is a topological isomorphism onto its closed range. Sometimes we additionally require homomorphisms to be contractive; otherwise we only assume boundedness, and we state this explicitly.

\begin{theorem}\label{thm:main}
There are no separable universal objects for the categories of separable Banach algebras, in either the commutative or the general setting, for contractive or merely bounded homomorphic embeddings. More precisely:
\begin{enumerate}
    \item[(A)] there is no separable \emph{commutative} Banach algebra $B$ such that every separable \emph{commutative} Banach algebra $A$ admits a \emph{contractive} homomorphic embedding $A\hookrightarrow B$;
    \item[(B)] there is no separable \emph{commutative} Banach algebra $B$ such that every separable \emph{commutative} Banach algebra $A$ admits a bounded homomorphic embedding $A\hookrightarrow B$;
    \item[(C)] there is no separable Banach algebra $B$ such that every separable Banach algebra $A$ admits a \emph{contractive} homomorphic embedding $A\hookrightarrow B$;
    \item[(D)] there is no separable Banach algebra $B$ such that every separable Banach algebra $A$ admits a bounded homomorphic embedding $A\hookrightarrow B$.
\end{enumerate}
\end{theorem}

\begin{remark}\label{rem:unital}
The theorem is stated in the standard non-unital category. A unital analogue follows by replacing $A_X$ with its unitisation: any unital embedding of the unitisation restricts to an injective homomorphism on $A_X$.
\end{remark}

\begin{remark}\label{rem:stronger}
As the proof will show, the closed-range hypothesis is not needed. In fact, there is no separable Banach algebra into which every separable Banach algebra admits even a bounded injective homomorphism; the same strengthening holds in the commutative and contractive settings.
\end{remark}

We present a single proof that covers all four cases. Contractivity plays no role beyond boundedness.

The key external input comes from the work of Johnson and Szankowski on complementably universal Banach spaces. Their original theorem shows that no separable Banach space is complementably universal for all separable Banach spaces \cite{JS1976}. For the present paper we use the following more precise real version from \cite[Theorem~2.3]{JS2009}.

\begin{theorem}[Johnson--Szankowski]\label{thm:JS}
Let
\[
\mathcal D_0^{\R}:=\{Y^*: Y\text{ is a closed subspace of }c_0(\R)\}.
\]
There is no separable real Banach space that is complementably universal for $\mathcal D_0^{\R}$.
\end{theorem}

The next lemma converts this real statement into the complex form needed below.

\begin{lemma}\label{lem:complex-reduction}
Let
\[
\mathcal D_0^{\C}:=\{X^*: X\text{ is a closed subspace of }c_0(\C)\}.
\]
If a complex Banach space $E$ is complementably universal for $\mathcal D_0^{\C}$, then the underlying real Banach space $E_{\R}$ is complementably universal for $\mathcal D_0^{\R}$.
\end{lemma}

\begin{proof}
Let $Y$ be a closed subspace of $c_0(\R)$, and let $Y_{\C}$ be its complexification. Then $Y_{\C}$ is a closed subspace of $c_0(\C)$. By assumption, $(Y_{\C})^*$ is isomorphic to a complemented subspace of $E$, and hence $((Y_{\C})^*)_{\R}$ is isomorphic to a complemented subspace of $E_{\R}$.

Now consider the bounded real-linear maps
\[
J:Y\to (Y_{\C})_{\R},\qquad J(y)=y+i0,
\]
and
\[
P:(Y_{\C})_{\R}\to Y,\qquad P(y_1+iy_2)=y_1.
\]
Then $P\circ J=\id_Y$, so $Y$ is complemented in $(Y_{\C})_{\R}$. Dualising, we obtain bounded operators
\[
P^*:Y^*\to ((Y_{\C})_{\R})^*,\qquad J^*:((Y_{\C})_{\R})^*\to Y^*,
\]
with $J^*P^*=\id_{Y^*}$. Therefore $Y^*$ is isomorphic to a complemented subspace of $((Y_{\C})_{\R})^*$.

Finally, $((Y_{\C})_{\R})^*$ is canonically isometrically isomorphic to $((Y_{\C})^*)_{\R}$ via
\[
\Phi:((Y_{\C})^*)_{\R}\to ((Y_{\C})_{\R})^*,\qquad \Phi(\varphi)=\operatorname{Re}\varphi.
\]
Hence $Y^*$ is isomorphic to a complemented subspace of $((Y_{\C})^*)_{\R}$, and therefore of $E_{\R}$. Since $Y\subseteq c_0(\R)$ was arbitrary, $E_{\R}$ is complementably universal for $\mathcal D_0^{\R}$.
\end{proof}

Recall that if $E$ and $U$ are Banach spaces and $\id_E$ factors through $U$, say
\[
\id_E = SR \qquad (R:E\to U,\ S:U\to E),
\]
then $R$ is bounded below, because $\norm{x}\le \norm{S}\,\norm{Rx}$ for every $x\in E$, and $RS$ is a bounded projection from $U$ onto $R(E)$. Conversely, if $J:E\to U$ is an isomorphic embedding and $P:U\to J(E)$ is a bounded projection, then
\[
\id_E = J^{-1}PJ.
\]
Thus factorisation of the identity on $E$ through $U$ is equivalent to complemented embeddability of $E$ in $U$.

Our argument is unified and uses a simple test-algebra construction. For each closed subspace $X$ of $c_0$ we build a separable commutative Banach algebra
\[
A_X := X^* \oplus_1 X \oplus_1 \C,
\]
whose multiplication records the canonical pairing between $X^*$ and $X$. Any injective homomorphism $\phi:A_X\to B$ yields bounded operators $p:X^*\to B$ and $S:B\to X^*$ with $S\circ p = \id_{X^*}$. Hence $X^*$ is isomorphic to a complemented subspace of $B$. If $B$ were universal for homomorphic embeddings of all separable Banach algebras, then the underlying complex Banach space of $B$ would be complementably universal for $\mathcal D_0^{\C}$, contradicting Lemma~\ref{lem:complex-reduction} and Theorem~\ref{thm:JS}.

\section{The test algebras}

Let $X$ be a closed subspace of $c_0$. Since $X$ is separable and $X^*$ is a quotient of $\ell_1$, both $X$ and $X^*$ are separable. Define
\[
A_X := X^* \oplus_1 X \oplus_1 \C,
\]
with multiplication
\begin{equation}\label{eq:product}
(f,x,\alpha)(g,y,\beta) := \bigl(0,0,f(y)+g(x)\bigr)
\qquad (f,g\in X^*,\ x,y\in X,\ \alpha,\beta\in \C).
\end{equation}

\begin{prop}\label{prop:AX}
For every closed subspace $X$ of $c_0$, the algebra $A_X$ is a separable commutative Banach algebra, its product is bilinear and contractive, and $A_X^3=\{0\}$.
\end{prop}

\begin{proof}
The formula \eqref{eq:product} defines a bilinear multiplication. Since $A_X$ is the $\ell_1$-sum of three Banach spaces, it is complete. Separability is clear. Also,
\[
|f(y)+g(x)| \le \norm{f}\,\norm{y}+\norm{g}\,\norm{x}
\le (\norm{f}+\norm{x}+|\alpha|)(\norm{g}+\norm{y}+|\beta|),
\]
so the multiplication is contractive.

The multiplication is commutative because the right-hand side of \eqref{eq:product} is symmetric in the two factors. Every product belongs to $\{0\}\oplus\{0\}\oplus\C$, and multiplying once more gives zero. Hence $A_X^3=\{0\}$, and in particular multiplication is associative.
\end{proof}

The next lemma is the key point. It converts an injective homomorphism of $A_X$ into a factorisation of the identity on $X^*$ through the target Banach algebra itself.

\begin{lemma}\label{lem:factor-id}
Let $X$ be a closed subspace of $c_0$, let $B$ be a Banach algebra, and let
\[
\phi:A_X\longrightarrow B
\]
be an injective homomorphism. Then the identity operator on $X^*$ factors through $B$.
\end{lemma}

\begin{proof}
Define bounded linear maps
\[
p:X^*\to B,\qquad p(f):=\phi(f,0,0),
\]
\[
q:X\to B,\qquad q(x):=\phi(0,x,0),
\]
and
\[
j:\C\to B,\qquad j(\alpha):=\phi(0,0,\alpha).
\]
Since $\phi$ is injective, $j$ is injective, and therefore $e:=j(1)\neq 0$.

For $f\in X^*$ and $x\in X$, using \eqref{eq:product} we get
\begin{equation}\label{eq:pxqx}
p(f)q(x)
= \phi(f,0,0)\,\phi(0,x,0)
= \phi\bigl((f,0,0)(0,x,0)\bigr)
= \phi(0,0,f(x))
= f(x)e.
\end{equation}

Choose $\lambda\in B^*$ such that $\lambda(e)=1$ by the Hahn--Banach theorem. Define
\[
S:B\to X^*,
\qquad
(Sb)(x):=\lambda\bigl(b\,q(x)\bigr)
\qquad (b\in B,\ x\in X).
\]
This is a bounded linear operator because
\[
|(Sb)(x)| \le \norm{\lambda}\,\norm{b}\,\norm{q(x)}
\le \norm{\lambda}\,\norm{q}\,\norm{b}\,\norm{x}.
\]
Now \eqref{eq:pxqx} implies that for every $f\in X^*$ and $x\in X$,
\[
\bigl(Sp(f)\bigr)(x)
= \lambda\bigl(p(f)q(x)\bigr)
= \lambda\bigl(f(x)e\bigr)
= f(x).
\]
Hence $Sp(f)=f$ for all $f\in X^*$, that is,
\[
S\circ p = \id_{X^*}.
\]
Thus $\id_{X^*}$ factors through $B$. In particular, $pS$ is a bounded projection from $B$ onto $p(X^*)$, and $p(X^*)$ is a complemented subspace of $B$ isomorphic to $X^*$.
\end{proof}

\section{Proof of the main theorem}

\begin{proof}[Proof of Theorem~\ref{thm:main}]
Assume that $B$ is a separable Banach algebra satisfying one of the four universality properties in Theorem~\ref{thm:main}. We treat all four cases simultaneously.

Let $X$ be an arbitrary closed subspace of $c_0$. By Proposition~\ref{prop:AX}, the algebra $A_X$ is separable, commutative, and has contractive multiplication. Therefore, in each of the cases \textup{(A)}--\textup{(D)}, the universality assumption yields a homomorphic embedding of $A_X$ into $B$ with the required norm control. By Lemma~\ref{lem:factor-id}, the identity operator on $X^*$ factors through $B$; equivalently, $B$ contains a complemented subspace isomorphic to $X^*$.

Since $X\subseteq c_0(\C)$ was arbitrary, the Banach space underlying $B$ is complementably universal for $\mathcal D_0^{\C}$. This contradicts Lemma~\ref{lem:complex-reduction} and Theorem~\ref{thm:JS}. The contradiction proves Theorem~\ref{thm:main}.
\end{proof}

\subsection*{Acknowledgements}
The author thanks the referee warmly for an exceptionally careful and generous report, for identifying the flaw in the earlier argument, and for a number of suggestions that substantially improved both the exposition and the final form of the paper. The present version, and in particular its shorter and more direct proof, owes a great deal to those comments.

Funding received from NCN Sonata-Bis 13 (2023/50/E/ST1/00067) is acknowledged with thanks.

\end{document}